\documentclass{amsart}

\usepackage{hyperref}

\newcommand{\lie}{\mathfrak{g}}
\newcommand{\kac}{\hat{\mathfrak{g}}}
\newcommand{\Ch}{\mathbb C[[h]]}
\newcommand{\Uha}{U_h(\kac)}

\newcommand{\B}{\mathcal B_h}
\newcommand{\Bha}{\B(\kac)}
\newcommand{\ba}[1]{\overline{#1}}
\newcommand{\End}{\text{End}}

\begin{document}

\title[Quantum affine reflection algebras]{Quantum affine reflection algebras\\
of type $d_n^{(1)}$ and reflection matrices}

\author{Gustav W. Delius \and Alan George}

\address{Department of Mathematics\\University of York\\York YO10
5DD\\United Kingdom}

\begin{abstract}
Quantum affine reflection algebras are coideal subalgebras of
quantum affine algebras that lead to trigonometric reflection
matrices (solutions of the boundary Yang-Baxter equation). In this
paper we use the quantum affine reflection algebras of type
$d_n^{(1)}$ to determine new $n$-parameter families of
non-diagonal reflection matrices. These matrices describe the
reflection of vector solitons off the boundary in $d_n^{(1)}$
affine Toda field theory. They can also be used to construct new
integrable vertex models and quantum spin chains with open
boundary conditions.
\end{abstract}

\maketitle

\section{Introduction}

Solutions to the reflection equation (also known as the boundary
Yang-Baxter equation) \cite{Cherednik:1984vs,Ghoshal:1994tm},
called reflection matrices, are required for the construction of
quantum integrable models with boundaries \cite{Sklyanin:1988yz}.
The reflection equation at the boundary is similar to the
Yang-Baxter equation in the bulk, the study of which led to the
axiomatization of quantum groups \cite{Drinfeld:1985rx}.
Intertwiners of these quantum groups are solutions of the
Yang-Baxter equation (called R-matrices) \cite{Jimbo:1986ua}.

An analogous method to find reflection matrices using intertwiners
of reflection equation algebras \cite{Sklyanin:1988yz} embedded as
coideal subalgebras in the corresponding Yangian or quantum affine
algebra has recently been developed
\cite{Delius:2001yi,Delius:2001qh}.  While there is a general
construction of these subalgebras in terms of the universal
R-matrix \cite{Delius:2001qh}, it is of practical importance  to
have a simple set of generators for them.  Such simple generators
have been found for coideal subalgebras of Yangians
\cite{Delius:2001yi} and for coideal subalgebras of quantum affine
algebras \cite{Delius:2001qh,Mezincescu:1998nw}, which we refer to
as ``quantum affine reflection algebras''.  With these expressions
it is easy to solve the intertwining condition, which has already
led to several new reflection matrices
\cite{Delius:2002ad,MacKay:2002at,Baseilhac:2002kf}.

The purpose of this letter is to demonstrate the practical utility
of these quantum affine reflection algebras by using them to
derive hitherto unknown reflection matrices corresponding to the
vector representation of $so(2n)$.  These reflection matrices are
needed, for example, to describe the reflection of solitons in
$d_n^{(1)}$ affine Toda theory off an integrable boundary
\cite{Delius:2001qh}.

\section{Quantum affine algebras\label{qa}}
In this section we summarize some essential facts about quantum
affine algebras. More details can be found  in, for example,
\cite{Chari:1994}.

Let $\kac$ be an affine Lie algebra with generalized Cartan matrix
$(a_{ij})_{i,j=0,\dots,n}$. The quantum affine algebra $U_h(\kac)$
is the unital associative algebra over $\Ch$ with generators
$x_i^+,\ x_i^-,\ h_i$, $i=0,\dots,n$ and relations
\begin{gather} \label{defrel}
[h_i,h_j] = 0,\ \ \ ~~
[h_i,x_j^\pm] = \pm a_{ij}\, x_j^\pm , \nonumber\\
[x_i^+,x_j^-] =  \delta_{ij}
\ \frac{q_i^{h_i} - q_i^{-h_i}}{q_i- q_i^{-1}},\\
\sum_{k=0}^{1-a_{ij}} (-1)^k \left[ \begin{array}{c} 1-a_{ij} \\ k
\end{array} \right]_{q_i} (x_i^\pm)^k x_j^\pm
(x_i^\pm)^{1-a_{ij}-k} = 0 \qquad i \not= j. \nonumber
\end{gather}
Here $\left[\begin{array}{c}a\\b\end{array}\right]_q$ are the
$q$-binomial coefficients. We have defined $q_i = e^{d_i h}$ where
the $d_i$ are such that $d_i a_{ij}$ is a symmetric matrix.

The
Hopf algebra structure of $\Uha$ is given by the
comultiplication $\Delta:\Uha \rightarrow \Uha \otimes \Uha$
defined by
\begin{align}
\Delta(h_i) &= h_i \otimes 1 + 1 \otimes h_i, \notag\\
\Delta(x_i^\pm) &= x_i^\pm \otimes q_i^{-h_i/2} +
 q_i^{h_i/2} \otimes x_i^\pm,
\end{align}
and the antipode $S$ and counit $\epsilon$ defined by
\begin{equation}
S(h_i)= - h_i, \ \ \ S(x_i^\pm) = - q_i^{\mp 1}\,x_i^\pm ,\ \ \
\epsilon(h_i) = \epsilon(x_i^\pm) = 0.
\end{equation}

For any $(n+1)$-tuple $\{\sigma_0,\dots,\sigma_n\}$ of invertible
constants in $\Ch$ there is a Hopf-algebra automorphism $\sigma$
of $\Uha$ defined by
\begin{equation}\label{resc}
    \sigma:\ x_i^\pm\mapsto\sigma_i^{\pm1}\,x^\pm_i,\ \ \
    h_i\mapsto h_i.
\end{equation}

Given some finite-dimensional $\Uha$-module $V$ with
representation map $\pi:\Uha\rightarrow\End(V)$ one can define a
spectral parameter dependent representation
$\pi_x:\Uha\rightarrow\End(V_x)$ by
\begin{equation}
    \pi_x=\pi\circ\sigma_x
\end{equation}
where $\sigma_x$ is an automorphism of the form \eqref{resc} with
the choice $\sigma_0=x$, $\sigma_i=1, i=1,\dots,n$. The tensor
product modules $V_{x_1}\otimes V_{x_2}\otimes\dots\otimes
V_{x_n}$ are irreducible for generic values for the spectral
parameters $x_i$.

Let $\check{R}(x)$ be the unique module homomorphism
$\check{R}(x):V_x\otimes V_1\mapsto V_1\otimes V_x$, i.e., let
$\check{R}(x)$ be a solution of the intertwining equation
\begin{equation}\label{ri}
    \check{R}(x)\,(\pi_x\otimes\pi_1)(\Delta(Q))=
    (\pi_1\otimes\pi_x)(\Delta(Q))\,\check{R}(x),~~~~
    \text{for all }~Q\in\Uha.
\end{equation}
Then $\check{R}(x)$ automatically satisfies \cite{Jimbo:1986ua}
the Yang-Baxter equation
\begin{equation}\label{ybe}
    (\check{R}(y)\otimes 1)(1\otimes\check{R}(xy))(\check{R}(x)\otimes
    1)=(1\otimes\check{R}(x))(\check{R}(xy)\otimes
    1)(1\otimes\check{R}(y)).
\end{equation}

\section{Quantum affine reflection algebras}

In this paper we are concerned with subalgebras $\Bha\subset\Uha$
generated by
\begin{equation}\label{qhat}
    \hat{Q}_j=q_j^{h_j/2}\left(x_j^+ + x_j^-\right) + \hat{\epsilon}_j
    \left(q_j^{h_j}-1\right),~~~~~j=0,\dots,n,
\end{equation}
for certain choices of the parameters $\hat{\epsilon}_j\in\Ch$. We
refer to these algebras as quantum affine reflection algebras. It
was shown in \cite{Delius:2001qh} how they are related to
Sklyanin's reflection equation algebras \cite{Sklyanin:1988yz}.
The algebras $\Bha$ are left coideal subalgebras of $\Uha$ because
\begin{equation}\label{qd}
    \Delta(\hat{Q}_j)=\hat{Q}_j\otimes 1+q_j^{h_j}\otimes\hat{Q}_j\in
    \Uha\otimes\Bha.
\end{equation}
Generically, the irreducible evaluation modules $V_x$ of $\Uha$
and their tensor products are also irreducible as modules of
$\Bha$.

A rescaling automorphism $\sigma$ of the form \eqref{resc} maps
the subalgebra $\Bha$ to a different subalgebra of $\Uha$ which is
isomorphic to $\Bha$ as an algebra and as a $\Uha$-coideal. In
particular the sign of $\hat{\epsilon}_j$ in \eqref{qhat} can be
changed by such an automorphism with $\sigma_j=-1$.

Let us assume that for some constant $\eta$ there exists an
intertwiner $K(x):V_{\eta x}\rightarrow V_{\eta/x}$, i.e., a
matrix $K(x)$ satisfying
\begin{equation}\label{ki}
    K(x)\,\pi_{\eta x}(\hat{Q})=\pi_{\eta/x}(\hat{Q})\,K(x),~~~~~~~~\forall\,\hat{Q}\in\Bha.
\end{equation}
Then $K(x)$ is automatically a solution of the reflection equation
\begin{equation}\label{bybe}
    (1\otimes K(y))\check{R}(xy)(1\otimes K(x))\check{R}(x/y)=
    \check{R}(x/y)(1\otimes K(x))\check{R}(xy)(1\otimes K(y))
\end{equation}
This is so because both sides of the equation are intertwiners
between the irreducible modules $V_{\eta x}\otimes V_{\eta y}$ and
$V_{\eta/x}\otimes V_{\eta/y}$. By Schur's lemma the two sides are
proportional. That they are equal then follows from the fact that
they have the same determinant and are equal at $h=0$.

\section{Vector representation of $d_n^{(1)}$}

We now specialize to the case of the $2n$-dimensional
representation of $\kac=d_n^{(1)}=\widehat{so(2n)}$. We choose
$d_i=1$ in \eqref{defrel}. The $\Uha$ generators are represented
by
\begin{align}
    \pi(x^+_j)&=E_{j,j+1}-E_{\ba{j+1},\ba{j}},&
    \pi(x^-_j)&=E_{j+1,j}-E_{\ba{j},\ba{j+1}},\\
    \pi(h_j)&=E_{j,j}-E_{j+1,j+1}-E_{\ba{j},\ba{j}}+E_{\ba{j+1},\ba{j+1}},\notag\\
    \pi(x^+_n)&=E_{n-1,\ba{n}}-E_{n,\ba{n-1}},&
    \pi(x^-_n)&=E_{\ba{n},n-1}-E_{\ba{n-1},n},\\
    \pi(h_n)&=E_{n,n}-E_{\ba{n},\ba{n}}-E_{\ba{n-1},\ba{n-1}}+E_{n-1,n-1},\notag\\
    \pi(x^+_0)&=E_{\ba{2},1}-E_{\ba{1},2},&
    \pi(x^-_0)&=E_{1,\ba{2}}-E_{2,\ba{1}},\\
    \pi(h_0)&=E_{\ba{2},\ba{2}}-E_{1,1}-E_{2,2}+E_{\ba{1},\ba{1}},\notag
\end{align}
where $E_{j,k}$ is the matrix with a $1$ in the $j$-th row and the
$k$-th column and all other entries $0$,
and
\begin{equation}
    \bar{j}=2n+1-j.
\end{equation}

The intertwining equation \eqref{ri} was solved by Jimbo
\cite{Jimbo:1986ua} to give the trigonometric R-matrix
\begin{align}
    \check{R}(x)=&(x-q^{-2})(x-\xi)\sum_{a\neq\ba{a}} E_{a,a}\otimes
    E_{a,a}+(x-1)(x-\xi)/q\sum_{a\neq b,\ba{b}} E_{a,b}\otimes
    E_{b,a}\notag\\
    &+(1-q^2)(x-\xi)\left(\sum_{a<b,a\neq\ba{b}}+x\sum_{a>b,a\neq\ba{b}}
    \right)E_{a,a}\otimes E_{b,b}\\
    &+\sum_{a\neq b} b_{ab}(x) E_{\ba{b},a}\otimes
    E_{b,\ba{a}},\notag
\end{align}
where $\xi=q^{2-2n}$ and
\begin{equation}
    b_{ab}(x)=\left\{
    \begin{array}{ll}
        (q^{-2}x-\xi)(x-1)&\text{for }~ a=b\\
        (q^{-2}-1)(\xi q^{b'-a'}(x-1)-\delta_{a,\ba{b}}(x-\xi))&\text{for }~ a<b\\
        (q^{-2}-1)x(q^{b'-a'}(x-1)-\delta_{a,\ba{b}}(x-\xi))&\text{for }~
        a>b
    \end{array}
    \right.,
\end{equation}
with
\begin{equation}\label{ip}
    a'=\left\{\begin{array}{ll}
        a+1/2&\text{if }~a\leq n\\
        a-1/2&\text{if }~a>n
    \end{array}\right.
\end{equation}

A Hopf-algebra automorphism $\sigma$ of the form \eqref{resc}
corresponds to a change of basis and a change of spectral
parameter, i.e.,
\begin{equation}
    \pi_x(\sigma(Q))=\Sigma^{-1}\,\pi_{\zeta x}(Q)\,\Sigma\ \ \
    \forall Q\in\Uha,
\end{equation}
where
\begin{align}
    &\zeta=\sigma_0\sigma_1(\sigma_2\dots\sigma_{n-2})^2\sigma_{n-1}\sigma_n,\notag\\
    &\Sigma=\text{diag}(\Sigma_1,\dots,\Sigma_{2n})\\
    &\Sigma_j=({\sigma_j}\dots{\sigma_{n-2}})^{-1},\qquad
    \Sigma_{\bar{j}}={\sigma_j}\dots{\sigma_{n}},\qquad
    j=1,\dots,n-2\notag\\
    &\Sigma_{n-1}=1,\qquad \Sigma_{n}=\sigma_{n-1},\qquad
    \Sigma_{\ba{n}}=\sigma_n,\qquad
    \Sigma_{\ba{n-1}}=\sigma_{n-1}\sigma_n.\notag
\end{align}
The R-matrix satisfies
\begin{equation}\label{rinv}
    (\Sigma^{-1}\otimes\Sigma^{-1})\check{R}(x)(\Sigma\otimes\Sigma)=\check{R}(x).
\end{equation}

\section{Reflection matrix for $d_n^{(1)}$}

It is now straightforward to solve the intertwining equation
$\eqref{ki}$ for the vector representation of $d_n^{(1)}$ to
obtain the solutions of the reflection equation \eqref{bybe}. One
finds that a solution exists only if
\begin{align}\label{ep}
    \eta&=\pm\,(-1)^n\,q^{1-n},&
    \hat{\epsilon}_j&=\pm i\, (q^{1/2}-q^{-1/2})^{-1},~~~~j=0,\dots,n.
\end{align}
The signs of the $\hat{\epsilon}_j$'s can be changed by a
rescaling automorphism, as discussed above, and the sign of $\eta$
can be changed by $x\rightarrow -x$. We thus now choose the upper
signs in \eqref{ep}. We write the solution as
\begin{equation}
    K(x)=\sum K_{ab}(x)\,E_{a,b}.
\end{equation}
We only need to give the entries on and above the antidiagonal
because the others are determined by the symmetry
\begin{equation}
    K_{\ba{a}\ba{b}}(x)=(-1)^n \widetilde{K}_{ab}(1/x),
\end{equation}
where the tilde indicates the change $q\mapsto 1/q$. The entries
are
\begin{equation}\label{k}
    K_{ab}(x)=\left\{\begin{array}{ll}
        i^{-\underline{a}-\underline{b}}\,q^{b'-n}\,k(x)&
        \text{for }~a+b\leq 2n \text{ and }a\neq b,\ba{b}\\
        i^{-\underline{a}-\underline{b}}\,q^{b'-n}\,(k(x)+(-1)^n\tilde{k}(1/x))(q+1)^{-1}&
        \text{for }~b=\ba{a}\text{~and~}a\leq n\\
        i^{-\underline{a}-\underline{b}}\,q^{b'-n}\,k(x)
        -(-1)^n\,q^{1/2}\,k(x)\\\qquad\qquad+
        ((-1)^n q^{(n-1)/2}+q^{(1-n)/2})&
        \text{for }~b=a
    \end{array}\right.
\end{equation}
where
\begin{equation}
    \underline{a}=\left\{\begin{array}{ll}a&\text{for }~a\leq n\\
        \ba{a}&\text{for }~a>n,\end{array}\right.
\end{equation}
the notation $a'$ was defined in \eqref{ip}, and
\begin{equation}
    k(x)=q^{n/2-1}(1-x).
\end{equation}
Using the invariance \eqref{rinv} of the R-matrix with respect to
the automorphisms \eqref{resc} we can write an $n$-parameter
family of solutions $K^+_\sigma$ of the reflection equation
\eqref{bybe}
\begin{equation}
    K^+_\sigma(x)=\Sigma^{-1}\,K(x)\,\Sigma.
\end{equation}
A second family $K^-_\sigma(x)$ of solutions is obtained by
choosing the opposite sign for $\eta$ in \eqref{ep}, leading to
$K^-_\sigma(x)=K^+_\sigma(-x)$. Furthermore, multiplying any of
these solutions by an arbitrary function of $x$ gives another
solution because the reflection equation is linear.

\section{Discussion}

The main new result of this letter is the hitherto unknown
trigonometric reflection matrices \eqref{k} for the vector
representation of $so(2n)$. It proves the practical utility of the
quantum affine reflection algebras with their simple generators
\eqref{qhat} and makes clear the possibility of generalization to
other representations and other algebras, which is the topic of
continuing research \cite{Mol,Nep}.

In order to use the reflection matrices \eqref{k} to describe
soliton reflection in $d_n^{(1)}$ affine Toda field theory they
need to be transformed to the principal gradation. This is
achieved by an automorphism $\sigma$ of the form \eqref{resc} with
$\sigma_0=1, \sigma_i=x, i=1,\dots,n$. Furthermore the
normalisation of the reflection matrix has to be fixed by the
requirements of unitarity and crossing symmetry, as was done for
$a_2^{(1)}$ in \cite{Gandenberger:1998eb}.

An equation very similar to the classical limit of our
intertwining condition \eqref{ki} appeared in
\cite{Bowcock:1995vp} in the context of constructing classically
integrable boundary conditions in affine Toda field theory.

Coideal subalgebras of non-affine quantized enveloping algebras
$U_h(\lie)$ have been studied \cite{Sug,Letzter:1999} in the
context of quantum symmetric spaces and $q$-orthogonal
polynomials. Expressions similar to \eqref{qhat} appear there.

The properties of quantum affine reflection algebras and their
representation theory deserve further study. For example a
knowledge of how irreducible representations of $\Uha$ branch into
irreducible representations of $\Bha$ will have physical
applications to the spectrum of boundary states in integrable
models.

\def\cprime{$'$}
\providecommand{\href}[2]{#2}\begingroup\raggedright\endgroup

\end{document}